\documentclass[11pt]{article}
\usepackage[cp850]{inputenc}
\usepackage{amssymb}
\usepackage{latexsym}
\def\Z{\mathbb{Z}}
\def\N{\mathbb{N}}
\def\R{\mathbb{R}}
\def\s{\frak{s}}

\begin{document}
\setlength{\parindent}{0pt}
\setlength{\parskip}{0.4cm}

\begin{center}

\large{\bf The Reciprocity Law for Dedekind Sums via the constant Ehrhart coefficient}
\footnote{This paper appeared in {\it American Mathematical Monthly} {\bf 106}, no. 5 (1999), 459-462. }

\normalsize{\sc Matthias Beck}

\end{center}



{\bf 1. Introduction}

The Dedekind sum can be defined for two relatively prime positive integers $a, b$ by
\[ \s (a,b) = \frac{ 1 }{ 4b } \sum_{ k=1 }^{ b-1 } \cot\frac{ \pi k a }{ b } \cot\frac{ \pi k }{ b } .  \]
These sums appear in various branches of mathematics: Number Theory, Algebraic Geometry, and Topology; 
they have consequently been studied extensively in various contexts. These include the quadratic reciprocity 
law (\cite{meyer}), random number generators (\cite{knuth}), group actions on complex manifolds (\cite{hirz}), 
and lattice point problems (\cite{pomm}, \cite{robins}).
Dedekind was the first to show the following reciprocity law (\cite{dedekind}):
\begin{equation}\label{recipr} 
  \s (a,b) + \s (b,a) = - \frac{ 1 }{ 4 } + \frac{ 1 }{ 12 } \left( \frac{ a }{ b } + \frac{ 1 }{ ab } + \frac{ b }{ a }  \right)  
\end{equation}
He was led naturally to this reciprocity law by considering the $\eta$-function 
$ \eta (\tau) = e^{ \frac{ \pi i \tau }{ 12 }  } \Pi_{ m=1 }^{ \infty } \left( 1 - e^{ 2 \pi i m \tau }  \right) $ 
on the complex upper half plane and transforming it under the action of the modular group $SL_{ 2 } ( \Z ) $.

Gau\ss 's law of quadratic reciprocity, for example, follows easily from (\ref{recipr}) (\cite{meyer}, \cite{grosswald}). 
We note that $ \s (a,b) = \s (a \mbox{ mod } b, b) $. Combining this with the reciprocity law (\ref{recipr}), 
one obtains a polynomial-time algorithm for computing $ \s (a,b) $, similar in spirit to the Euclidean algorithm. 
From this point of view, it is not surprising (though not obvious) that $ \s (a,b) $ can be 
expressed efficiently in terms of the continued fraction expansion of $a/b$ ({\cite{hick}, \cite{vardi}).

Rademacher was one of the pioneers in the use of Dedekind sums (\cite{white}); in fact, he found several proofs 
of (\ref{recipr}) (\cite{grosswald}). We present yet another proof, which establishes a simple connection 
with lattice point enumeration in polytopes.  The reciprocity law (\ref{recipr}) follows readily once the reader is familiar with the computation of the 
coefficients of the Ehrhart polynomial for a lattice polytope. 

\vspace{1cm}

{\bf 2. Counting Lattice Points}

Let $\Z^{ n } \subset \R^{ n } $ be the $n$-dimensional integer lattice, and $\cal P$ an $n$-dimensional 
lattice polytope in $ \R^{ n } $. So $\cal P$ is a compact simplicial complex of pure dimension $n$ whose 
vertices lie on the lattice. For $t \in \N$, denote by $ L ( {\cal P} , t ) $ the number of lattice points 
in the closure of the dilated polytope $ t {\cal P} := \{ tx : x \in {\cal P} \} $. Ehrhart (\cite{ehrhart}) proved that 
$ L ( {\cal P} , t ) $ is a polynomial in $t$ of degree $n$. Moreover, 
\[ L ( {\cal P} , t ) = \mbox{Vol}({\cal P}) t^{ n } + \frac{ 1 }{ 2 } \mbox{Vol}(\partial {\cal P}) t^{ n-1 } + \dots + \chi ( {\cal P}) . \]
Here, Vol($\partial {\cal P}$) denotes the surface area of ${\cal P}$ normalized with respect to 
the sublattice on each face of ${\cal P}$, and $\chi ( {\cal P}) $ is the Euler characteristic 
of ${\cal P}$. We note that, for convex polytopes $\cal P$, $ \chi ( {\cal P}) = 1 $ (\cite{ehrhart}). 

In this paper, we focus on the case $\R^{2}$, where Ehrhart's result is known as Pick's Theorem (\cite{grun}, \cite{diaz}):
For a convex lattice polytope ${\cal P} \in \R^{2}$, 
\[ L ( {\cal P} , t ) = A t^{ 2 } + \frac{1}{2} B t + 1 , \]
where $A$ is the area and $B$ the number of boundary lattice points of ${\cal P}$. 

In the general case, the other coefficients of $ L ( {\cal P} , t ) $ are not as easily accessible. In 
fact, until quite recently a method of computing these coefficients was unknown.
There has been recent progress in this direction (\cite{barv}, \cite{brion}, \cite{kantor}, \cite{puk}); 
Diaz and Robins (\cite{robins}) found a way of proving a cotangent representation for the 
generating function $ \sum_{ t=0 }^{ \infty } L ( {\cal P} , t ) e^{ -2 \pi s t } $, thereby 
deriving a formula for the Ehrhart coefficients of $ L ( {\cal P} , t ) $.
For our purposes, the following result (a straightforward consequence of Corollary 1 in \cite{robins}) is sufficient: 

{\bf Theorem.} {\it
  Let ${\cal P}$ denote the simplex in $\R^{ n } $ with the vertices \\
  $ (0, \dots , 0), \left( a_{ 1 } , 0 , \dots 0 \right) , \left( 0, a_{ 2 } , 0, \dots , 0 \right) , \dots , \left( 0, \dots , 0, a_{ n }  \right) $, 
  where $ a_{ 1 } , \dots , a_{ n } \in \N $ are pairwise coprime. Denote the corresponding 
  Ehrhart polynomial by $ L ( {\cal P} , t ) = \sum_{ j=0 }^{ n } c_{ j } t^{ j } $. Then 
  $c_{ m }$ is the coefficient of $ s^{ -(m+1) } $ in the Laurent expansion at $s=0$ of 
  \begin{eqnarray*} \frac{ \pi^{ m+1 }  }{ m! 2^{ n-m } p } 
    \sum_{ r=1 }^{ p } \left( 1 + \coth \frac{ \pi }{ a_{ 1 }  } ( s + i r ) \right) \left( 1 + \coth \frac{ \pi }{ a_{ 2 }  } ( s + i r ) \right) \cdot \dots \\
         \cdot \left( 1 + \coth \frac{ \pi }{ a_{ n }  } ( s + i r ) \right) \left( 1 + \coth \frac{ \pi }{ p } ( s + i r ) \right)  , \end{eqnarray*}
  where $ p = a_{ 1 } \cdot \dots \cdot a_{ n } $. }

Note that the appearance of cotangent products in this result leads us to expect Dedekind sums 
in some form within the coefficients of the Ehrhart polynomial, thus also within the formulas 
for the number of lattice points in simplices. In fact, the nontrivial cases of dimension three 
(\cite{rad}) and four (\cite{rosen}) involve classical Dedekind sums. Both formulas can easily 
be obtained through the above Theorem. 

We will use this result in an indirect way. Precisely, we will compute $c_{0}$ according to the 
above theorem, and make use of the fact that $c_{ 0 } = \chi ( {\cal P}) = 1$.
Dedekind's reciprocity law (\ref{recipr}) will follow from this idea if we consider the case of dimension $n=2$.

\vspace{1cm}


{\bf 3. Proof of the Reciprocity Law}

According to the theorem, we have to find the coefficient of $s^{ -1 }$ of the Laurent series at $s=0$ of (note that $(a,b)=1$)
\begin{equation}\label{sum}
  \frac{ \pi }{ 4ab } \sum_{ r=1 }^{ ab } \left( 1 + \coth \frac{ \pi }{ a } ( s + i r ) \right) \left( 1 + \coth \frac{ \pi }{ b } ( s + i r ) \right) \left( 1 + \coth \frac{ \pi }{ ab } ( s + i r ) \right)
\end{equation}
The Laurent expansion of each factor depends on $r$:
\[ 1 + \coth \frac{ \pi }{ c } ( s + i r ) = \left\{ \begin{array}{ll}
   S_{ c } := \frac{ c }{ \pi } s^{ -1 } + 1 + \frac{ \pi }{ 3c } s + O \left( s^{ 3 }  \right) & \mbox{ if } c | r \\
   R_{ c } := 1 + \coth \frac{ \pi i r }{ c } + O ( s ) & \mbox{ if } c \not| r \end{array} \right. \]
To keep track of the various cases, we introduce the notation 
\[ \chi_{ c } = \left\{ \begin{array}{l}
   1 \mbox{ if } c | r \\
   0 \mbox{ if } c \not| r , \end{array} \right. \]
so that we can write $ 1 + \coth \frac{ \pi }{ c } ( s + i r ) = S_{ c } \chi_{ c } + R_{ c } \left( 1 - \chi_{ c }  \right) $, 
and (\ref{sum}) becomes
\[ \sum_{ r=1 }^{ ab } \left( S_{ a } \chi_{ a } + R_{ a } \left( 1 - \chi_{ a }  \right) \right) \left( S_{ b } \chi_{ b } + R_{ b } \left( 1 - \chi_{ b }  \right) \right) \left( S_{ ab } \chi_{ ab } + R_{ ab } \left( 1 - \chi_{ ab }  \right) \right) . \]

Now, expand this into all 8 terms, and consider each summand according to the number of $S_{ c }$ factors:
\begin{enumerate} 
  \item Terms with one $S_{ c }$ factor are 
    \begin{eqnarray}\label{S_a}
          S_{ a } \chi_{ a } R_{ b } \left( 1 - \chi_{ b }  \right) R_{ ab } \left( 1 - \chi_{ ab }  \right) 
          &=& S_{ a } R_{ b } R_{ ab } \chi_{ a } \left( 1 - \chi_{ b } - \chi_{ ab } + \chi_{ ab }  \right) \nonumber \\
          &=& S_{ a } R_{ b } R_{ ab } \left( \chi_{ a } - \chi_{ ab } \right) \end{eqnarray}
    and, similarly,
    \begin{equation}\label{S_b}
          R_{ a } \left( 1 - \chi_{ a }  \right) S_{ b } \chi_{ b } R_{ ab } \left( 1 - \chi_{ ab }  \right) 
          = R_{ a } S_{ b } R_{ ab } \left( \chi_{ b } - \chi_{ ab } \right) . \end{equation}
    The summand with $S_{ ab }$ is zero 
    (note that $\chi_{ a } \chi_{ ab } = \chi_{ b } \chi_{ ab } = \chi_{ ab } $, and $\chi_{ a } \chi_{ b } = \chi_{ ab } $). 
    To compute the contribution of (\ref{S_a}), note that the support of 
    $ \chi_{ a } - \chi_{ ab } $ in $ \left\{ 1, \dots , ab \right\} $ is 
    $ \left\{ ka : 1 \leq k \leq b - 1 \right\} $; thus its contribution to (\ref{sum}) is
    \begin{eqnarray*}
      \frac{ \pi }{ 4 ab } &\cdot& \frac{ a }{ \pi } \ \sum_{ k=1 }^{ b-1 } \left( 1 + \coth \frac{ \pi i ka }{ b } \right) \left( 1 + \coth \frac{ \pi i ka }{ ab } \right) \\
      &=& \frac{ 1 }{ 4b } \sum_{ k=1 }^{ b-1 } \left( 1 - i \cot \frac{ \pi ka }{ b } \right) \left( 1 - i \cot \frac{ \pi k }{ b } \right) \\
      &=& \frac{ 1 }{ 4b } \sum_{ k=1 }^{ b-1 } 1 - \cot \frac{ \pi ka }{ b } \cot \frac{ \pi k }{ b } + i \dots \\
      &=& \frac{ 1 }{ 4 } - \frac{ 1 }{ 4b }  - \s (a,b) .
    \end{eqnarray*}
    The imaginary part in the preceding sum has to be zero, because the original generating function is real.
    Similarly, (\ref{S_b}) gives a contribution of $ \frac{ 1 }{ 4 } - \frac{ 1 }{ 4a }  - \s (b,a) $.
  \item There are no terms with two $S_{ c }$ factors, because 
    \[ S_{ a } \chi_{ a } S_{ b } \chi_{ b } R_{ ab } \left( 1 - \chi_{ ab }  \right) = S_{ a } S_{ b } R_{ ab } \chi_{ ab } \left( 1 - \chi_{ ab }  \right) = 0 \]
    and
    \[ S_{ a } \chi_{ a } R_{ b } \left( 1 - \chi_{ b }  \right) S_{ ab } \chi_{ ab } = S_{ a } R_{ b } S_{ ab } \chi_{ ab } \left( 1 - \chi_{ b }  \right) = 0 . \]
  \item Finally, the term $ S_{ a } \chi_{ a }  S_{ b } \chi_{ b } S_{ ab } \chi_{ ab } = S_{ a } S_{ b } S_{ ab } \chi_{ ab } $ 
    has support $ \{ ab \} $, and gives a contribution of 
    \begin{eqnarray*}
      \frac{ \pi }{ 4ab } \left( \frac{ a }{ \pi } \frac{ b }{ \pi } \frac{ \pi }{ 3ab } + \frac{ a }{ \pi } \frac{ ab }{ \pi } \frac{ \pi }{ 3b } + \frac{ b }{ \pi } \frac{ ab }{ \pi } \frac{ \pi }{ 3a } +  \frac{ a }{ \pi } + \frac{ b }{ \pi } + \frac{ 
ab }{ \pi } \right) \\
      = \frac{ 1 }{ 12 } \left( \frac{ 1 }{ ab } + \frac{ a }{ b } + \frac{ b }{ a }  \right) + \frac{ 1 }{ 4 } \left( \frac{ 1 }{ b } + \frac{ 1 }{ a } + 1 \right) .
    \end{eqnarray*}
\end{enumerate}
Adding all contributions, we arrive at 
\[ 1 = c_{ 0 } = \frac{ 3 }{ 4 } + \frac{ 1 }{ 12 } \left( \frac{ 1 }{ ab } + \frac{ a }{ b } + \frac{ b }{ a }  \right) - \s (a,b) - \s (b,a) , \]
the desired reciprocity law (\ref{recipr}).

The same method applied to dimension $n=3$ does not give any further results. However, for $n=4$, 
higher dimensional Dedekind sums (\cite{zagier}) appear within the computations, so that this case 
is likely to provide new results. I will discuss this matter in a future paper. 

\vspace{1cm}

{\bf Acknowledgements}

The author would like to thank Sinai Robins for helpful suggestions and invaluable support.

\vspace{1cm}

\footnotesize
\nocite{*}
\addcontentsline{toc}{subsubsection}{References}
\bibliography{thesis}
\bibliographystyle{alpha}

\vspace{1cm}

\sc Dept. of Mathematics\\
Temple University\\
Philadelphia, P.A. 19122 

{\tt matthias@math.temple.edu}\\
{\tt http://www.math.temple.edu/$\sim$matthias}

\end{document}